\documentclass[11pt]{article}
\usepackage{amssymb,amsmath,latexsym, bbm}

\newtheorem{thm}{Theorem}[section]
\newtheorem{cor}[thm]{Corollary}
\newtheorem{lemma}[thm]{Lemma}
\newtheorem{prop}[thm]{Proposition}
\newtheorem{defn}[thm]{Definition}

\newtheorem{remark}[thm]{Remark}
\numberwithin{equation}{section}

\def\E {{\mathbb E}}
\def\P {{\mathbb P}}
\def\R {{\mathbb R}}

\def\N {{\mathbb N}}

\def\LL{\mathcal L}
\def\1{{\mathbf 1}}

\def\<{{\langle}}
\def\>{{\rangle}}
\def\eps{\varepsilon}
\def\wh{\widehat}

\def\proof{{\medskip\noindent {\bf Proof. }}}

\def\qed{{\hfill $\square$ \bigskip}}

\begin{document}

\allowdisplaybreaks
\bibliographystyle{plain}

\title{\bf One-dimensional heat equation with discontinuous conductance}

\author{{\bf Zhen-Qing Chen}\thanks{Research partially supported
by NSF Grant DMS-1206276 and  NNSFC Grant 11128101.} \quad  and \quad {\bf Mounir Zili}}

\maketitle

\begin{abstract}

We study a second-order parabolic equation  with divergence form elliptic operator, having piecewise constant diffusion coefficients with two points of discontinuity. Such partial differential equations  appear in the modelization of diffusion  phenomena in medium consisting of three kind of materials.
Using probabilistic methods, we present an explicit expression of the fundamental solution under certain conditions.  We also derive small-time asymptotic expansion of the PDE's solutions in the general case. The obtained results are  directly usable in applications.

\end{abstract}

\bigskip

\noindent {\bf AMS 2000 Mathematics Subject Classification}: Primary
 60H10, 60J55; Secondary 35K10, 35K08

\bigskip\noindent
{\bf Keywords and phrases}:  Stochastic differential equations, semimartingale local time, strong solution, heat kernel,  asymptotic expansion

\bigskip

\section{Introduction}

In many practical applications, one is encountered with heat propagations
in heterogeneous media consisting of different kind of materials.
Mathematically, such kind of heat propagation is modeled by heat equations
with divergence form elliptic operators having
discontinuous diffusion coefficients. In this paper we exam one-dimensional model where the medium consisting of three different kind of materials.

Let $a$, $a_i, \rho_i$, $i=1, 2, 3$ be positive constants and define
$$A(x)=  a_1 \1_{\{ x \le 0 \} } + a_2 {\bf 1}_{\{ 0 < x \le a \} } + a_3 {\bf 1}_{\{ a < x  \} } \quad \hbox{and} \quad
\rho(x)= \rho_1 \1_{\{ x \le 0 \} } + \rho_2 {\bf 1}_{\{ 0 < x \le a \} }
+ \rho_3 {\bf 1}_{\{ a < x  \} }.
$$
Set
\begin{equation} \label{e:1.1}
 \LL = \frac{1}{2\rho (x)} \frac{d}{dx} \left( \rho (x)A(x) \frac{d}{dx} \right),
\end{equation}
which is a self-adjoint operator in $L^2(\R; \rho (x)dx)$.
The parabolic equation
\begin{equation}\label{e:1.2}
 \frac{\partial u(t, x)}{\partial t} = \LL u(t, x)
\end{equation}
describes the heat propagation in the line that consists of three different kind of material.
It is well-known that there is an $m$-symmetric diffusion process $X$
associated with $\LL$, where $m(dx):= \rho (x) dx$.
Since both $A(x)$ and $\rho (x)$ are bounded between two positive
constants, it is known (see \cite{GSC}) that $X$ possesses a jointly continuous
transition density function $q(t, x, y)$ with respect to the measure $m(dx)$
on $\R$ that is symmetric in $x$ and $y$ and enjoys the following
Aronson type Gaussian estimates: there exists positive constants $C_1, C_2 \geq 1$
so that for every $t>0$ and $x, y\in \R$,
\begin{equation}\label{e:1.3}
C_1^{-1}  t^{-1/2} \exp (-C_2 |x-y|^2/t) \leq q(t, x, y) \leq
C_1 t^{-1/2} \exp (- |x-y|^2/(C_2t)) .
\end{equation}
Using Dirichlet form theory,
 the following semimartingale representation of $X$ is derived
 in \cite[Proposition 5]{LejMart}:
 \begin{equation}\label{e:1.4}
dX_t = \sqrt{A(X_t)} dB_t +\frac{\rho_2 a_2- \rho_1 a_1}{\rho_2 a_2+ \rho_1 a_1}
\, d \wh L^0_t + \frac{\rho_3 a_3- \rho_2 a_2}{\rho_3 a_3+ \rho_2 a_2}
\, d \wh L^a_t ,
\end{equation}
where $B$ is one-dimensional Brownian motion,
$\wh L^0$ and $\wh L^a$ are symmetric semimartingale local times of $X$ at $0$ and $a$, respectively. Here for a  semimartingale $Y$, its
symmetric semimartingale local time  $\wh L^w_t(Y)$ at  $w\in \R$ is
 defined to be
\begin{equation}  \label{e:10}
 \wh L^w_t (Y):= \lim_{\eps\to 0} \frac1{2 \eps} \int_0^t {\bf 1}_{[w-\eps, w+\eps]} (Y_s) \, d\<Y\>_s,
 \end{equation}
 where $\< Y\>$ is the predictable quadratic variation process of $Y$.
 We define the semimartingale local time for $Y$ at level $w\in \R$ to be
\begin{equation}  \label{e:1.5}
 L^w_t (Y)=\lim_{\eps\to 0}
\frac1{ \eps} \int_0^t {\bf 1}_{[w, w + \eps]} (Y_s)\, d\<Y\>_s.
 \end{equation}
It is known (see \cite{BC}) that equation \eqref{e:1.4} admits a unique strong solution. Moreover, $X$ is a strong solution to
\eqref{e:1.4} if and only if it is a strong solution to
\begin{equation}\label{e:1.6}
dX_t = \sqrt{A(X_t)} dB_t +\frac12 \left( 1- \frac{\rho_1 a_1}{\rho_2 a_2}
\right)  \, d  L^0_t (X) + \frac12 \left( 1-\frac{\rho_2 a_2}{\rho_3 a_3} \right)
\, d L^a_t (X) .
\end{equation}

The purpose of this paper is to study the parabolic equation \eqref{e:1.2}
with initial value $u(0, x)= h(x)$ and solve it as explicitly as possible,
where $h$ is a piecewise $C^{N+1}$ function on $\R$ of the form
\begin{equation}\label{e:1.8}
h(x)= h_1(x)\1_{\{x\leq 0\}}+ h_2(x) \1_{\{ 0<x\leq a\}}+ h_3 (x)\1_{\{ x>a\}},
\end{equation}
with $h_i$ being bounded with continuously differentiable up to order  $N+1$,
$i=1, 2, 3$, for some integer $N \geq 1$.
 This kind of PDEs appears in the modelization of diffusion phenomena in many areas, for example, in   geophysics \cite{Lej}, ecology \cite{Can},   biology \cite{Nic}
 and so on. The non-smoothness of the coefficients reflects the multilayered medium.
The solution $u(t, x)$ of \eqref{e:1.2} with initial value $h(x)$ admits a probabilistic representation:
\begin{equation}\label{e:1.9}
 u(t, x)=\E_x [ h(X_t)], \quad t\geq 0.
\end{equation}
Here the subscript $x$ in $\E_x$ means the expectation is taken with respect to
the law of the diffusion process $X$ that starts from $x$ at time $t=0$.
In the first part of this paper, we show that, when either
$\rho_1 \sqrt{a_1}=\rho_2 \sqrt{a_2}$ or $\rho_2 \sqrt{a_2}=\rho_3 \sqrt{a_3}$,
we can reduce the diffusion process of \eqref{e:1.3} to
skew Brownian motion through transformations by scale functions.
Thus in this case, we are able to derive an explicit formula for
the transition density function $q(t, x, y)$  of the diffusion process $X$ with respect to the measure $m(dx)$. The kernel $q(t, x, y)$ is
the fundamental solution (also called heat kernel) $q(t, x, y)$ of $\LL$.
Skew Brownian motions are a subclass of diffusion processes $X$ of \eqref{e:1.3}
with $a_1=a_2=a_3=1$ and $\rho_2=\rho_3$. It is first introduced by Ito and McKean
in 1963, and has since been studied extensively by many authors; see, for example,
\cite{ABTWW, BC, RY, WA} and the references therein.
However we do not know the explicit formula for $q(t, x, y)$
in the general case.
In \cite{GAV}, formulas for the fundamental solutions are derived
in terms of inverse Fourier transform and Green functions, but they
are not explicit.

Though we do not have explicit formula for its heat kernel,
in the second part of this paper, we are able to
derive asymptotic expansion in small time of
solutions of \eqref{e:1.2} in the  general case
with piecewise $C^{N+1}$ initial data,
by making use of the explicit heat kernel obtained
under the special case $(a_2, \rho_2)=(a_3, \rho_3)$.
The small-time asymptotic expansion we get in this paper is explicit and can be directly used in applications.
This extends the work of \cite{Zi} and of \cite{Zil}, where only one discontinuity
is studied. In  \cite{ETOR} and \cite{LejMart}, numerical methods are proposed to study \eqref{e:1.4} with discontinuous coefficients.

The rest of the paper is organized as follows.
In the section 2, we derive  explicit formula for the fundamental solution
 of $\LL$ in the case when either $\rho_1 \sqrt{a_1}=\rho_2 \sqrt{a_2}$
 or $\rho_2 \sqrt{a_2}=\rho_3 \sqrt{a_3}$.

 We then study the small-time asymptotic expansion of the solution of \eqref{e:1.2}
 in the general case for any positive constants $a$, $a_i$ and $\rho_i$.
Extensions to the more general piecewise constant coefficient case with more
than two discontinuities are mentioned at the end of the paper.

\section{Fundamental solution}

In view of the probabilistic representation \eqref{e:1.9}, properties of
solution to the heat equation \eqref{e:1.9} can be deduced from
the properties of the diffusion process $X$ of \eqref{e:1.6}.
For notational simplicity, let's rewrite SDE \eqref{e:1.6}
as
\begin{equation}\label{e:2.1}
 \begin{cases}
 dY_t^x = \Big( p {\bf 1}_{\{ Y_t^x \le 0 \} } + q {\bf 1}_{\{ 0 < Y_t^x \le a \} } + r {\bf 1}_{\{ a < Y_t^x  \} } \Big)
dB_t + \frac{\alpha}2 d  L_t^0 (Y^x) + \frac{\beta}2  d   L_t^a (Y^x) ,\\
Y^x_0 = x \in \R,
\end{cases}
 \end{equation}
 where $p, q, r$ are positive constants,  $\alpha, \beta \in (-\infty, 1)$ and $B$  a one-dimensional Brownian motion. The superscript $x$ indicates  the process
 $Y^x$ starts from $x$ at $t=0$.
In the particular case when $p = q = r =1$ and $\beta = 0$, $Y_t^x$ is just a Skew Brownian motion of parameter $(2-\alpha )^{-1}$ (cf. \cite{WA}), and in the case when $p = q = r =1$, $Y_t^x$ is the double-skewed Brownian motion (cf. \cite{RA}).

The next two theorems are known; see  \cite[Theorems 2.1 and 2.2]{BC}
or \cite{LG}.
We give a proof here not only for reader's convenience but also
some formulas in the proof will be used later to derive the transition density
function of $Y$.

\begin{thm}\label{T:1}
 For every $\alpha < 1 $ and $\beta < 1 $,
SDE \eqref{e:2.1}  has a unique strong solution for every $x \in \R $.
\end{thm}

\proof   Define
\begin{equation} \label{e:2.2}
s(x)=\begin{cases}
x/p &\hbox{for } x< 0, \\
x/q &\hbox{for } x \in [0, a], \\
(x-a)/r + a/q &\hbox{for } x > a.
\end{cases}
\end{equation}
Clearly, $s$ is strictly increasing and one-to-one.
Let $\sigma$ denote the inverse of $s$:
\begin{equation} \label{e:2.3}
\sigma(x)=\begin{cases}
px &\hbox{for } x<0, \\
qx &\hbox{for } x \in [0, s(a)], \\
r(x-s(a)) + s(a)q &\hbox{for } x > s(a).
\end{cases}
\end{equation}
Let $s'_\ell$ and $\sigma'_\ell$ denote  the
left hand derivative of $s$ and $\sigma $, respectively; that is,
 $$s'_\ell(x)=\begin{cases}
1/p &\hbox{for } x\leq 0, \\
1/q &\hbox{for } x \in (0, a], \\
1/r &\hbox{for } x > a,
\end{cases}
\quad \hbox{and} \quad
\sigma'_\ell (x)=\begin{cases}
p &\hbox{for } x\leq 0, \\
q &\hbox{for } x \in (0, s(a)], \\
r &\hbox{for } x > s(a),
\end{cases}
$$
The second derivative  $\sigma''(dx) = (q-p)\delta_{\{0\}}(dx)
+(r-q)\delta_{\{ s(a) \}} (dx)$. Here $\delta_{\{z\}}$ denotes
the Dirac measure concentrated at point $z$.

Since both $\frac{q(\alpha -1)}p+1 $ and $\frac{r(\beta-1)}{q}+1$
are strictly less than $1$, by Theorem $2.1$
of \cite{BC}, the following SDE has a unique strong solution
 $Z=\{ Z_t, \, t\geq 0 \}$ with $Z_0 =s(x)$:

\begin{equation}\label{e:6}
dZ_t=  dB_t + \frac12 \left(\frac{q(\alpha -1)}p+1\right) dL^0_t(Z)+\frac12
 \left(\frac{r(\beta-1)}{q}+1\right) dL^{s(a)}_t (Z) .
\end{equation}

  Let $Y= \sigma(Z)$.
We  will show that $Y$ is a solution to \eqref{e:2.1}. Clearly
$Y_0 = x$.  Function $\sigma$ is a piecewise linear function and can
be expressed as a difference of two convex functions.
By Ito-Tanaka's formula (see \cite[Theorem VI.1.5]{RY}), $Y$
is a semimartingale and in fact, for $t\geq 0$,
\begin{eqnarray} \label{e:7}
Y_t&=& \sigma(Z_0) + \int_0^t \sigma'_\ell(Z_s) dZ_s
+ \frac12 \int L^w_t(Z) \sigma''(dw) \nonumber\\
&=& x +  \int_0^t \Big( p {\bf 1}_{\{ Y_t \le 0 \} } + q {\bf 1}_{\{ 0 < Y_t \le a \} }
 + r {\bf 1}_{\{   Y_t>a  \} } \Big) dB_t
 \nonumber\\
&&  + \frac{\alpha\, q}2 L^0_t(Z)+\frac{\beta\, r}2 L^{s(a)}_t(Z).
 \end{eqnarray}
By  Corollary VI.1.9 of Revuz and Yor \cite{RY}, for $t \geq 0$,
\begin{eqnarray} \label{e:8}
L^{s(a)}_t(Z)&=&\lim_{\eps\downarrow 0} \frac{1}{\eps} \int_0^t 1_{[s(a),s(a)+\eps)}(Z_s) d\<Z\>_s
  = \lim_{\eps\downarrow 0} \frac{1}{\eps} \int_0^t 1_{[a,a+r\eps)}(Y_s) ds
  \nonumber \\
 &=&  \frac1{r} \lim_{\eps\downarrow 0} \frac{1}{r \eps} \int_0^t 1_{[a,a+r\eps)}(Y_s) d\<Y\>_s
 = L^a_t(Y)/r.
 \end{eqnarray}
A similar calculation shows that $L^0_t(Z)= L^0_t(Y)/q$.
Thus we have from \eqref{e:7} that $Y= \sigma(Z)$ is a strong solution
for  \eqref{e:2.1} with $Y_0= x$.

We now examine pathwise uniqueness. Suppose that
$Y^x$ is a solution of \eqref{e:2.1}.
Since $s$ is the difference of two convex functions, we can
apply Ito-Tanaka's formula to get
\begin{eqnarray} \label{e:9}
s(Y_t^x) = s(x)+\int_0^t s'_\ell (Y_s^x)   dY_s^x +
  \frac12 \int_{\R} L^w_t(Y^x) s''(dw).
\end{eqnarray}
Let $Z= s(Y^x)$. A similar calculation as above shows that
$Z$ satisfies SDE \eqref{e:6} with initial value $Z_0= s(x)$.
Since by Theorem $2.1$ of \cite{BC} solutions to \eqref{e:6} is unique,
hence so is solution $Y^x$ to \eqref{e:2.1}.
\qed

\begin{thm} \label{T:2.2}
The process $Y^x$ is a strong solution of
 the SDE \eqref{e:2.1},  if and only if $Y^x$ is a strong solution to
\begin{equation}\label{e:11}
   dY_t^x = \Big( p {\bf 1}_{\{ Y_t^x \le 0 \} } + q {\bf 1}_{\{ 0 < Y_t^x \le a \} } + r {\bf 1}_{\{ a < Y_t^x  \} } \Big)
dB_t  +   \frac{\alpha}{2-\alpha}   d \wh L_t^0 (Y^x) +
\frac{\beta }{2-\beta}   d   \wh L_t^a (Y^x)
 \end{equation}
 with $Y_0^x=x$.
\end{thm}

 \proof This is due to the fact that, for solution $Y^x$ of \eqref{e:2.1} (see  the proof of Theorem 2.2 in \cite{BC}),
\begin{equation} \label{e:12}
  L^0_t(Y^x)= \frac{2}{2-\alpha} \wh L^0_t(Y^x),
 \qquad L^a_t(Y^x)= \frac{2}{2-\beta} \wh L^a_t(Y^x),
\end{equation}
and $ L^w_t(Y^x)=\wh L^w(Y^x)$ for  $w\notin \{0, a\}$.
\qed

\begin{cor} \label{C:3}
For every $p, q > 0, \alpha<1 $ and $x \in \R $, the following SDE
  has a unique strong solution for every $x \in \R$:
  \begin{equation}\label{e:13}
 \begin{cases}  dX_t^{x} = \Big( p {\bf 1}_{\{ X_t^{x} \le 0 \} } +
 q {\bf 1}_{\{ 0 < X^{x}_t\}}  \Big)
dB_t + \frac{\alpha}2 d  L_t^0 (X^{x})  ,\\
X_0^x = x \in \R .
\end{cases}
 \end{equation}
Moreover, the transition probability density of the diffusion $X^{x}$ is
given by
\begin{equation} \label{e:14}
 \begin{array}{rcl}
\displaystyle  p^X(t, x, y)
 &= & \displaystyle  \frac{1}{\sqrt{2 \pi t}}
 \left( \frac{ {\bf 1}_{ \{ y \le 0 \} } }{p}
+ \frac{ {\bf 1}_{ \{ y > 0 \} } }{q} \right) \times \left\{ \exp \left( - \frac{(f(x) - f(y))^2}{2 t}\right) \right. \\
\noalign{\vskip 2mm}
& & + \displaystyle \left.
\frac{p+q (\alpha -1)}{p-q (\alpha -1)}\, {\rm sign} (y)\, \exp \left( - \frac{(\mid f(x) \mid  + \mid f(y) \mid)^2}{2t}\right) \right\}
\end{array}
\end{equation}
where $f(y)= \frac{y}{p}\1_{\{y\leq 0\}}+\frac{y}{q} \1_{\{y>0\}}$.
\end{cor}

\proof Taking  $q = r$ and $\beta = 0$, we have by Theorem \ref{T:1}   that
for every $p, q > 0$, $ \alpha < 1$ and $ x \in \R $,
 \begin{equation}\label{e:15}
 dX^{x}_t = \Big( p {\bf 1}_{\{ X^{x}_t \le 0 \} } + q {\bf 1}_{\{ 0 < X^{x}_t  \} }  \Big) dB_t + \frac{\alpha}2 dL_t^0 (X^{x})  ,
\end{equation}
has a unique strong solution with $X^{x}_0={x}$.
Let $s$ be the function defined by \eqref{e:2.2}.
We see from \eqref{e:6} that
$Z:= s(X^{x})$ satisfies SDE \eqref{e:6} with $q=r$ and $\beta = 0$; that is,
\begin{equation}\label{e:16}
dZ_t= dB_t  + \frac12 \left(\frac{q(\alpha -1)}p+1\right) dL^0_t(Z).
\end{equation}
Thus $Z$ is a skew Brownian, whose transition density function is explicitly known. Using Theorem \ref{T:2.2}, we can rewrite $Z$ in terms of its symmetric semimartingale local time $\wh L^0(Z)$:
\begin{equation} \label{e:17}
 dZ_t= dB_t  +  \frac{p + q(\alpha -1)}{p-q (\alpha -1)} \,   d \wh L^0_t(Z).
 \end{equation}
Thus by  \cite{WA}  (see also Exercise III.(1.16) on page 82 of \cite{RY}),
$Z$ has a transition density function $p^Z(t, x, y)$ with respect to the Lebesgue measure on $\R$ given by
\begin{eqnarray*}
  p^Z(t, x,y) &=& \frac{1}{\sqrt{2 \pi t}}
  \left( \exp \left( - \frac{(x - y)^2}{2 t}\right)+
\frac{p+q (\alpha -1)}{p-q (\alpha -1)}\, {\rm sign} (y)\, \exp \left( - \frac{(|x| + |y|)^2}{2t}\right) \right).
\end{eqnarray*}
The desired formula \eqref{e:14} now follows from the fact that
$$
p^X(t, x,y) = \frac{1}{\sigma _l'(s(y))} p^Z(t, s(x),s(y))
$$
  \qed

\begin{remark}\label{R:2.4}  \rm \begin{description}
\item{(i)} When $\alpha=1-(p^2/q^2)$, Corollary \ref{C:3} in particular
recovers the formula  given in \cite[Theorem 5.1]{ABTWW}
(see the proof of Proposition 5.1 of \cite{ABTWW} for a more
explicit one)
for the transition density function of $X^x$ in \eqref{e:13}.
Note that in this case, the more compact formula \eqref{e:14}
in this paper
is equivalent to the form given in \cite{ABTWW}.
Using the connection to skew Brownian motion
established in the proof of Corollary \ref{C:3},
one can derive the joint distribution of
$X^x$ given by \eqref{e:11}, its semimartingale (respectively, symmetric semimartingale)
local time $L^0(X^x)$ (respectively, $\wh L^0(X^x)$) at level $0$ and its
occupation time on $[0, \infty)$ from that of the corresponding
skew Brownian motion. The latter can be found in \cite[Theorem 1.2]{ABTWW}.
Same remark applies to the diffusion process $Y^x$
satisfying SDE \eqref{e:21} below.

\item{(ii)} The diffusion process $X$ of \eqref{e:13} can be called {\it
oscillating skew Brownian motion}, or
{\it skewed oscillating Brownian motion} as suggested by Wellner \cite{We}.
When $ \alpha =1- (p^2/q^2)$, the marginal distribution of $X^x$ of \eqref{e:13}
(characterized by its transition density function $p^X(t, x, y)$)
is the Fechner distribution  in statistics; see \cite{We}.
This corresponds exactly to the case when the density function
$p^X(t, 0, y)$ is continuous in $y$.
For generator $\LL$ in \eqref{e:1.1} or its associated SDE \eqref{e:1.6},
the condition $ \alpha =1- (p^2/q^2)$ corresponds precisely to the case
of $\rho_2=\rho_3$.  \qed
\end{description}
\end{remark}

The next theorem allows us to obtain
 the transition density function for solution $Y$ of \eqref{e:2.1}
 when $\beta= 1-\frac{q}{r}$.

 \begin{thm}\label{T:4}
 Let $p, q, r, a>0$ and $\alpha<1$.
 For $x\in \R$, let $X^x$ denote the unique strong solution of \eqref{e:13}
 with initial value $X^x_0=x$.
 The process $Y$ defined by:
\begin{equation} \label{e:20}
  Y_t =    \begin{cases}
  a + \frac{r}{q} (X^x_t - a) ^+ - ( X_t^x - a)^- & \hbox{if }   Y_0 = x \le a \\
  \noalign{\vskip 3mm}
  a + \frac{r}{q} (X^{\frac{q}{r}(x-a) + a}_t - a) ^+ - ( X_t^{\frac{q}{r}(x-a) + a} - a)^- & \hbox{if }   Y_0 = x > a .
\end{cases}
 \end{equation}
 is the unique strong solution of the stochastic differential equation:
 \begin{equation}\label{e:21}
dY^x_t  = \Big( p {\bf 1}_{\{ Y_t^x \le 0 \} } + q {\bf 1}_{\{ 0 < Y_t^x \le a \} } + r {\bf 1}_{\{ a < Y_t^x  \} } \Big)
dB_t + \frac{\alpha}2 d  L_t^0 (Y^x) + \frac12 \left(1-\frac{q}{r}\right)
d   L_t^a (Y^x)
\end{equation}
with initial value $Y^x_0=x$.
 \end{thm}

\proof We will present the proof just in the case where $\displaystyle x \le a$. The other case is similar.
Consider the bijective function defined by:
$$  \varphi  (x) = a+ (r/q) (x-a)^+ - (x-a)^- \qquad \hbox{for } x\in \R.$$
Then  $ Y  = \varphi (X^x)$. Note that $\varphi (x)=x$ for $x\leq a$,
and $\varphi (x)>a$ if and only if $x>a$.
Since $  \varphi $ is the difference of two convex functions, applying the Ito-Tanaka formula to $Y=\varphi (X^x)$, we obtain:
\begin{eqnarray}  \label{e:22}
 Y_t  &=& Y_0+ \int_0^t \varphi'_\ell (X^x_s) dX^x_s
  +\frac12 \int_{\R} L^w_t(X^x) \varphi '' (dx) \nonumber \\
  &=& x+ \int_0^t \Big(  {\bf 1}_{\{ X_s^x \leq a \} } + \frac{r}{q}
{\bf 1}_{\{ X_s^x > a \} } \Big) dX^x_s + \frac12 \left(\frac{r}{q}-1\right)
L^a_t(X^x) \nonumber \\
&=& x+\int_0^t \Big( p {\bf 1}_{\{ X^x_s \le 0 \} } +
 q {\bf 1}_{\{ 0 < X^x_s \le a \} } + r {\bf 1}_{\{ X^x_s>a  \} } \Big)
dB_t \nonumber \\
&& +  \frac{\alpha}2 d  L_t^0 (X^x) +  \frac12 \left(\frac{r}{q}-1\right)
L^a_t(X^x) \nonumber \\
&=& x+\int_0^t \Big( p {\bf 1}_{\{ Y^x_s \le 0 \} } +
 q {\bf 1}_{\{ 0 < Y^x_s \le a \} } + r {\bf 1}_{\{   Y^x_s>a  \} } \Big)
dB_t \nonumber  \\
&&  + \frac{\alpha}2 d  L_t^0 (X^x) + \frac12 \left(\frac{r}{q}-1\right)
L^a_t(X^x).
\end{eqnarray}
By the same reasoning as that for (\ref{e:8}), we have
 \begin{equation} \label{e:23}
 L^0_t(X^x)=L^0_t(Y) \quad \hbox{and} \quad  L^a_t(X^x)=\frac{q}{r}L^a_t(Y).
\end{equation}
This together with \eqref{e:22} shows that $Y$ is a strong solution to SDE \eqref{e:21} with $Y_0=x$. The uniqueness follows from Theorem \ref{T:1}.
\qed

\begin{cor} \label{C:5}
When $\beta= 1-\frac{q}{r}$, the solution $Y$ to SDE \eqref{e:2.1} has a
  transition probability density $p^Y(t, x, y)$ with respect to the Lebesgue
  measure on $\R$ given by
\begin{equation} \label{e:24}
 \begin{array}{rcl}
\displaystyle  p_t^Y(x,y) &= & \displaystyle  \frac{1}{\sqrt{2 \pi t}} \Bigg( \frac{ {\bf 1}_{ \{ y \le 0 \} }}{p}
+ \frac{ {\bf 1}_{ \{ 0 < y \le a \} }}{q} + \frac{ {\bf 1}_{ \{ y>a \} } }{r} \Bigg) \times  \Bigg\{ \exp \Bigg( - \frac{(s(x) - s(y))^2}{2 t}\Bigg) \\
\noalign{\vskip 2mm}
&&+ \, \displaystyle
\frac{p+q (\alpha -1)}{p-q (\alpha -1)}\, {\rm sign} (y)\, \exp \Bigg( - \frac{(\mid s(x) \mid  + \mid s(y) \mid)^2}{2t}\Bigg) \Bigg\}
\end{array}
\end{equation}
\end{cor}

\proof This follows immediately from Theorem \ref{T:4} and Corollary \ref{C:3}. \qed

\begin{remark} \rm
SDE \eqref{e:2.1} with $\beta= 1-\frac{q}{r}$  corresponds exactly to
SDE \eqref{e:1.6} with $\rho_2 \sqrt{a_2} =\rho_3 \sqrt{a_3}$.
If we use $q(t, x, y)$ to denote the transition density function of $Y$
with respect to its symmetrizing measure $m(dx)=\rho (x)dx$, then clearly
we have
$$ q(t, x, y)= p^Y(t, x, y)/ \rho (y).
$$
\end{remark}

\section{Small-time asymptotic expansion}

  In this section, we will give an explicit asymptotic expansion, in small time,
of the function $\displaystyle  u(x,t) = \E(h(Y_t^x))$, in the general case where
$p,q,r>0$ and $ \alpha , \beta   \in (- \infty , 1)$.

\begin{defn}
For every $ k \in \N$, we denote by $\displaystyle {\bf erfc_k}$ the function defined on $\R$ by:
$$   {\bf erfc_k}(z) = \frac{2}{\sqrt{\pi } } \int_z^{+ \infty } u^k e^{-u^2}du,
\qquad z\in \R.
$$
\end{defn}

\medskip

Observe that ${\bf erfc_k}(0) = \frac{1}{\sqrt{ \pi }} \Gamma \Big( \frac{k+1}{2} \Big)$,  where $\Gamma $ is the Euler Gamma function. The
following elementary properties are straightforward from its  definition.

\begin{lemma} \label{L:8} For every given integer $n\geq 1$,

\begin{description}
\item{\rm (i)}
As $ z \rightarrow + \infty $,
$\displaystyle {\bf erfc_k} (z) = o (z^{-n}) $.

\item{\rm (ii)} As $z \rightarrow - \infty $:\hspace{2mm}
${\bf erfc_k} (z) = \frac{1}{\sqrt \pi } \Big( 1 + (-1)^k \Big) \Gamma \Big( \frac{k+1}{2} \Big) + o (|z|^{-n}) $.
\end{description}
\end{lemma}

Here the notation $\displaystyle o(|z|^{-n})$ as $\displaystyle z\to +\infty$ (respectively, $\displaystyle z\to - \infty$) means that $\displaystyle \lim_{z\to +\infty}
\frac{o(|z|^{-n})}{|z|^{-n}}=0$ (respectively,
$\displaystyle \lim_{z\to -\infty} \frac{o(|z|^{-n})}{|z|^{-n}}=0$).
Similar notation will also be used
for $o(t^n)$ as $t\to 0$.
Note that by Corollary \ref{C:3},  SDE
  \begin{equation}\label{e:55}
 \begin{cases}  dS_t^x = \Big( q {\bf 1}_{\{ S_t^x \le a \} } +
 r {\bf 1}_{\{ a < S^x_t \}}  \Big)
dB_t + \frac{\beta }2 d  L_t^a (S^x)   \\
S_0^x = x \in \R
\end{cases}
 \end{equation}
has a unique strong solution for every $x \in \R$ and
it has a transition probability density
with respect to the Lebesgue measure on $\R$ given by
\begin{equation} \label{e:56}
 \begin{array}{rcl}
\displaystyle  p (t, x,y) &= & \displaystyle  \frac{1}{\sqrt{2 \pi t}} \Bigg( \frac{ {\bf 1}_{ \{ y \le a \} } }{q}
+ \frac{ {\bf 1}_{ \{ y > a \} } }{r} \Bigg) \times  \Bigg\{ \exp \Bigg( - \frac{(g(x) - g(y))^2}{2 t}\Bigg) \\
\noalign{\vskip 2mm}
&  & + \, \displaystyle
\frac{q+r (\beta -1)}{q-r (\beta -1)}\, {\rm sign} (y-a)\, \exp \Bigg( - \frac{( |g(x)|  + | g(y)|)^2}{2t}\Bigg) \Bigg\},
\end{array}
\end{equation}
where $g(x)= \frac{x}{q} \1_{\{x \le a\}} +
\frac{x}{r} \1_{\{ x > a\}}$.

\bigskip

The following estimation is a key in small-time asymptotic expansion for the solution
of heat equation \eqref{e:1.2}. It reduces the case to diffusions of the form
\eqref{e:13} for which we have explicit knowledge about their heat kernels.

\begin{lemma} \label{L:15}
Define
 \begin{equation} \label{e:60}
  \Xi^{x} = \begin{cases}
 X^{x}, \hspace{1mm} solution \hspace{1mm} of \hspace{1mm} sde \hspace{1mm} (\ref{e:13})  & \hbox{if } \, \displaystyle x \le \frac{a}{2} \\
 \noalign{\vskip 3mm}
 S^{x}, \hspace{1mm} solution \hspace{1mm} of \hspace{1mm} sde \hspace{1mm} (\ref{e:55})  &\hbox{if } \, \displaystyle x > \frac{a}{2}.
 \end{cases}
 \end{equation}
Then there exist   positive constants $c_1$ and $ c_2$ such that for every $t > 0$,
 $$
 {\Bbb P} \Big(  \Xi_s \neq Y^x_s   \hbox{ for some } s\in [0, t]\Big)
  \le c_1 \exp\Big( -c_2  a^2/ t \Big) ,
 $$
where $Y^x$ is the diffusion process of \eqref{e:2.1}.
\end{lemma}

\proof  We will prove the lemma in the case where $\displaystyle x \le \frac{a}{2}$. The other case is similar.
For $t > 0$, let us denote by
$$
 \tau _t = \inf \{ t \ge 0; \Xi_s = a \} = \inf \{ t \ge 0; X_s^{x} = a \} ,
$$
and
\begin{equation} \label{e:63}
A_t = \Big\{ \omega\in \Omega:  \Xi_s(\omega) \neq Y^x_s(\omega)
 \hbox{ for some } s\in [0, t]\Big\} = \Big\{ \omega\in \Omega:    X_s^{x}(\omega) \neq Y^x_s(\omega)  \hbox{ for some } s\in [0, t] \Big\}.
\end{equation}
It is clear that, $\displaystyle A_t \subset \{ \tau _t \le t \}$ and consequently,
$$
  \P (A_t)   \le  \P (\tau _t \le t ) \leq
  \P_{x}\left(\sup_{0 \le s \le t} |X_s^{x} - x | \ge
  |a - x| \right)\leq \P_{x}\left(\sup_{0 \le s \le t} |X_s^{x} - x | \ge
  a/2 \right) .
$$
Here $\P_x$ is the probability law of $X^x$.
Since $X$ is a diffusion whose transition density function $q(t, x, y)$ with respect to the measure $m(dx)=\rho (x) dx$ enjoys the Aronson type Gaussian estimate \eqref{e:1.3}, we have by a similar argument as that for \cite[Lemma II.1.2]{St} that
$$ \P_{x}\left(\sup_{0 \le s \le t} |X_s^{x} - x | \ge
  a/2 \right) \leq c_1 e^{-c_2 a^2/ t}.
  $$
This proves the lemma. \qed

\begin{cor} \label{C:16}
Let  $\Xi^{x} $ be  defined by  \eqref{e:60}. Then
  for every integer $n\geq 1$ and time  $t>0$ small,
$$
 \E\Big[ h(Y_t^x)\Big]  = \E\Big[ h(\Xi^{x}_t) \Big] +
 o(t^n) .
$$
 \end{cor}

 \proof    Let us  consider the case where $\displaystyle x \le \frac{a}{2}$; the other case is similar.
For all $\displaystyle t \ge 0$,
$$
\E\Big[ h(Y_t^x)\Big] = \E\Big[ h(Y_t^x); A_t \Big]
+ \E\Big[ h(Y_t^x);  A_t^c \Big] ,
$$
where $A_t$ is the set defined in (\ref{e:63}), and $A^c$ is the complementary of $A$.
Since $X_s^x(w) = Y_s^x(w)$ for all $s\in [0, t]$ on $A_t^c$, we have
 \begin{eqnarray*}
 \E\Big[ h(Y_t^x)\Big] &=&   \E\Big[ h(X_t^x); A_t^c \Big] +
\E\Big[ h(Y_t^x); A_t  \Big] \\
&=&   \E\Big[ h(X_t^x) \Big]  +
\E\Big[ h(Y_t^x)-  h(X_t^x);  A_t  \Big]  .
\end{eqnarray*}
The desired result now follows from Cauchy-Schwarz inequality, the fact that
$h$ is bounded and Lemma \ref{L:15}.
 \qed

Let $h$ be a piecewise $C^{N+1}$  function of the form
\eqref{e:1.8}.
By Taylor expansion, for every $x \in \R \setminus \{ 0, a \}$,
$$h(y)= \sum_{j=1}^N \frac{h^{(j)}(x)}{j!} (y-x)^j + O(|y-x|^{N+1}),$$
and for $x\in \{ 0, a \}$ and $i=1, 2, 3$,
$$h_i(y)= \sum_{j=1}^N \frac{h_i^{(j)}(x)}{j!} (y-x)^j + O(|y-x|^{N+1}).$$
Here the notation $O(|y-x|^{N+1})$ means that
there are  constants $C,  \delta>0$ so that
the term $O(|y-x|^{N+1})$ is no larger than $C|y-x|^{N+1}$ for any
$y\in \R$ with $|y-x|<\delta$.
Thus we have the following.

\begin{prop} \label{P:3.5}
If $\displaystyle x \in {\Bbb R } \setminus \{ 0, a \}$,
\begin{equation} \label{eq:34}
\E \left[h(\Xi^{x}_t)\right] = \sum_{j=0}^N \frac{1}{j!}  h^{(j)}(x) \E \Big[(\Xi^{x}_t -x) ^j \Big] +  O \Big( |\Xi^{x}_t-x|^{N+1} \Big)  ,
\end{equation}
and if $x \in \{ 0, a \}$,
\begin{eqnarray}
\label{eq:35}
  \E [h(\Xi^{x}_t)] &=&   \sum_{j=0}^N \frac{1}{j!}  h^{(j)}_1(x) \E \Big[(\Xi^{x}_t-x)^j
   {\bf 1}_{\{ \Xi^{x}_t \le 0 \}} \Big]  +
   \sum_{j=0}^N \frac{1}{j!} h^{(j)}_2(x) \E \Big[(\Xi^{x}_t-x)^j
   {\bf 1}_{\{ 0< \Xi^{x}_t \le a\}} \Big] \nonumber \\
& +&   \sum_{j=0}^N \frac{1}{j!} h^{(j)}_3(x) \E \Big[  (\Xi^{x}_t-x) ^j {\bf 1}_{\{ \Xi^{x}_t > a \}} \Big]  +   O \Big( |\Xi^{x}_t-x|^{N+1} \Big)   .
\end{eqnarray}
\end{prop}

 Using  the transition probability densities of the diffusions $X^{x}$  and $S^x $ given in corollary \ref{C:3} and in \eqref{e:56}, we can compute
 the expectations in \eqref{eq:34}-\eqref{eq:35}.

\begin{prop} \label{P:17}
Let $k\geq 1$ be an integer.
\begin{description}
\item{\rm (i)} For $ x < 0$,
$$
 \E \left[ \left( \Xi_t^x -x  \right)^k  \right] =  (-1)^k p^k 2^{k/2-1} t^{k/2}  {\bf erfc_k}\Big( \frac{x}{p \sqrt{2t}} \Big)
 +     \sum_{j=0}^k  {k\choose j}  2^{j/2-1}t^{j/2} A_j (x)\, {\bf erfc_j}\Big( \frac{-x}{p\sqrt{2t}}\Big),
$$
where
$$
A_j (x) =  \frac{x^{k-j}}{p-q (\alpha -1)} \Big(    (p+q(\alpha -1)) p^{j} (-1)^{k+1}2^{k-j} + 2pq^j   (\frac{q}{p}-1)^{k-j}   \Big)
\quad \hbox{and} \quad  {k\choose j}= \frac{k!}{j! (k-j)!}.
$$

\item{\rm (ii)} For $ 0< x < a$,
$$
\E \left[ \left( \Xi_t^x -x  \right)^k  \right] =
2^{k/2-1}q^k t^{k/2}
{\bf erfc_k} \Big( \frac{-x}{q \sqrt{2t}} \Big)  +  \sum_{j=0}^k {k \choose j} 2^{j/2 -1}t^{j/2} B_j (x)\,{\bf erfc_j} \Big( \frac{x}{q \sqrt{2t}} \Big) ,
$$
where
$$
B_j (x) = \frac{x^{k-j}}{p-q(\alpha -1)} \Big( - 2q (\alpha -1) (-p)^j \Big( \frac{p}{q} -1 \Big) ^{k-j}  +  (p+ q(\alpha -1)) q^j (-2)^{k-j} \Big) .
$$

\item{\rm (iii)} For $ x > a$,
$$
 \E \left[ \left( \Xi_t^x -x  \right)^k  \right] =
2^{k/2-1} r^k t^{k/2} {\bf erfc_k} \Big( \frac{a-x}{r \sqrt{2t}} \Big)
 +   \sum_{j=0}^k {k \choose j} 2^{j/2 -1}t^{j/2} C_j (x) \, {\bf erfc_j} \Big( \frac{x-a}{r \sqrt{2t}} \Big) ,
$$
where
$$
C_j (x)= \frac{(x-a)^{k-j}}{q-r( \beta -1)} \Big( -2r(\beta -1)(-1)^jq^j \left(\frac{q}{r} -1 \right)^{k-j}  + \left(q+r(\beta -1)\right)r^j 2^{k-j} (-1)^{k-j} \Big).
$$

\end{description}

\end{prop}

As a consequence of the above Proposition, and by Lemma \ref{L:8}, we obtain the following expansion
of $\E \left[ \left( \Xi_t^x -x  \right) ^k  \right]$ for $x \in {\Bbb R} \setminus \{ 0, a \}$ and small time $t>0$.

\begin{cor} \label{C:18}
For any positive integers $n>k\geq 1$, $\displaystyle x \in {\Bbb R} \setminus \{ 0, a \} $, as $\displaystyle t \rightarrow 0^+$,
$$
 \E \left[ \left( \Xi_t^x -x  \right) ^k  \right] = 2^{k/2-1} D_k (x) t^{k/2}+
 o(t^n),
$$
 where
\begin{equation}\label{e:3.12}
D_k(x) = \frac{1}{\sqrt{\pi }}(1+(-1)^k) \Gamma \left( \frac{k+1}{2} \right) \left[ p^k {\bf 1}_{\{x < 0 \} } + q^k {\bf 1}_{ \{ 0 < x < a \} } + r^k {\bf 1}_{\{x > a \} } \right] .
\end{equation}
\end{cor}

Using again the transition probability densities of the diffusions $X^{x}$  and $S^x $,  we get

\begin{prop} \label{P:19}
For $x \in \{ 0, a \}$ and integer $k\geq 1$,
 \begin{eqnarray*}
 \E \left[ \left( \Xi_t^x -x \right)^k {\bf 1}_{\{ \Xi_t^x \le 0\}}  \right]
 & =&
 \begin{cases} \displaystyle \frac{q(1-\alpha )}{p-q(\alpha -1)} (-1)^k p^k \frac{(2t)^{k/2}}{\sqrt{\pi }} \Gamma \Big( \frac{k+1}{2} \Big)    &\hbox{when }
   x = 0  \\
 \displaystyle \frac{r(1- \beta )}{q-r(\beta  -1)} (-1)^k q^k (2t)^{k/2} {\bf erfc_k} \Big( \frac{a}{q \sqrt{2t}} \Big)  &\hbox{when } x = a,
\end{cases}  \\
  \E \left[ \left( \Xi_t^x -x \right)^k {\bf 1}_{\{0< \Xi_t^x \le a\}}\right]
  & =&
{\cal A}_k (x)\,  q^k (2t)^{k/2} \left( \frac{1}{\sqrt{\pi }} \Gamma \Big( \frac{k+1}{2} \Big) -   {\bf erfc_k} \Big( \frac{a}{q \sqrt{2t}} \Big) \right) , \\
 \E \left[ \left( \Xi_t^x -x \right)^k {\bf 1}_{\{\Xi_t^x > a\}} \right]
 &=&  \begin{cases}
  \displaystyle \frac{p}{p-q(\alpha -1)}  q^k 2^{k/2}t^{k/2} {\bf erfc_k} \Big( \frac{a}{q \sqrt{2t}} \Big)     &\hbox{ when }  x = 0  \\
 \displaystyle   \frac{q}{q-r(\beta  -1)} r^k \frac{2^{k/2}t^{k/2}}{\sqrt{\pi }} \Gamma \Big( \frac{k+1}{2} \Big)  &\hbox{ when } x = a .
  \end{cases}
\end{eqnarray*}
Here
\begin{equation}\label{e:3.19}
{\cal A}_k (x) =
\begin{cases}
 \displaystyle  \frac{p}{p-q(\alpha -1)}    &\hbox{ when }  x = 0 \\
\noalign{\vskip 3mm}
  \displaystyle \frac{r(1-\beta )(-1)^k }{q-r(\beta -1)}    &\hbox{ when } x =  a .
\end{cases}
\end{equation}
\end{prop}

From Proposition \ref{P:19} and Lemma \ref{L:8} we deduce:

\begin{cor} \label{C:20}
For $ x \in \{ 0, a \}$ and integers $n>k\geq 1$,
 \begin{eqnarray*}
\E \left[ \left( \Xi_t^x -x \right)^k {\bf 1}_{ \{\Xi_t^x \le 0\}} \right]
&=&  \begin{cases}
 \displaystyle \frac{q(1-\alpha )}{p-q(\alpha -1)} (-1)^k p^k \frac{2^{k/2}t^{k/2}}{\sqrt{\pi }} \Gamma \Big( \frac{k+1}{2} \Big)    &\hbox{ when }   x = 0  \\
\displaystyle  o(t^n) &\hbox{ when }  x = a;
\end{cases} \\
\E \left[ \left( \Xi_t^x -x \right)^k {\bf 1}_{ 0< \Xi_t^x \le a} \right]
 & =&
\displaystyle {\cal A}_k(x) q^k (2t)^{k/2}  \frac{1}{\sqrt{\pi }} \Gamma \Big( \frac{k+1}{2} \Big) + o(t^n)  , \\
\displaystyle \E \left[ \left( \Xi_t^x -x \right)^k {\bf 1}_{ \Xi_t^x > a} \right]
&=&   \begin{cases}
 \displaystyle  o(t^n)    & \hbox{ when }  x = 0  \\
\displaystyle  \frac{q}{q-r(\beta  -1)} r^k \frac{2^{k/2}t^{k/2}}{\sqrt{\pi }} \Gamma \Big( \frac{k+1}{2} \Big)  &\hbox{ when } x = a ,
\end{cases}
\end{eqnarray*}
where ${\cal A}_k (x)$ is defined in \eqref{e:3.19}.
\end{cor}

Combining Proposition \ref{P:3.5} with
Corollaries \ref{C:18} and \ref{C:20}, we have

 \begin{thm} \label{T:21}
 For small time and $x\in \R$,
$$
   \E \left[ h(\Xi_t^x) \right]=   \sum_{k=0}^N b_k(x) t^{k/2} + O (t^{(N+1)/2}),
$$
 where
\begin{equation}\label{e:87}
b_k(x) = \begin{cases}
\displaystyle 2^{k/2-1} D_k(x)   \frac{1}{k!} \frac{\partial ^k h(x)}{\partial x^k} &\hbox{ when }
 x \in {\Bbb R} \setminus \{ 0, a \}  \\
 \displaystyle  \frac{2^{k/2}}{k! (p-q(\alpha - 1)) \sqrt{\pi} } \Gamma ( \frac{k+1}{2})
\left(\frac{\partial ^k h_1(0)}{\partial x^k} q(1-\alpha )(-1)^kp^k +
\frac{\partial ^k h_2(0)}{\partial x^k} pq^k \right)
&\hbox{ when } x=0  \\
 \displaystyle \frac{2^{k/2}}{k! (q- r (\beta -1)) \sqrt{\pi} } \Gamma  ( \frac{k+1}{2})
\left(  \frac{\partial ^k h_2(a)}{\partial x^k} r(1- \beta ) (-1)^kq^k +
\frac{\partial ^i h_3(a)}{\partial x^k} qr^k \right)
&\hbox{ when } x=a
\end{cases}
\end{equation}
with the function $D_k(x)$ given by \eqref{e:3.12}.
 \end{thm}

By Theorem \ref{T:21} and Corollary \ref{C:16} we get:

 \begin{cor} \label{C:22} For every $x\in \R$, when $t>0$ is small,
$$
  \E [h(Y_t^x) ]=  \sum_{k=0}^N b_k(x) t^{k/2} + O (t^{(N+1)/2}),
$$
where the function $b_k (x)$ is defined by \eqref{e:87}.
\end{cor}

\begin{remark}\rm (i)  When $a = 0$ and $r = q$, we recover
  the expansion obtained in \cite{Zi}.

(ii) Employing the same approach   used in this section,
 one can obtain without any difficulty a similar small-time expansion
 for the solution of the heat equation \eqref{e:1.2} where
  $A(x)$ and $\rho (x)$ are piecewise constant with
  $n\geq 3$ points of discontinuity; that is,
  in case of the following SDE:
 $$
 \begin{cases}
 \displaystyle dY_t^x =  \sqrt{A(Y_t^x)}
dB_t + \sum_{i=1}^n \frac{\alpha _i }2 d  L_t^{a_i} (Y^x) ,\\
\displaystyle Y_0 = x \in \R,
\end{cases}
$$
  where $\displaystyle a_1=0 < a_2 < ... < a_n $,  $\displaystyle  p_i \in (0, + \infty )$, $\displaystyle  \alpha _i \in ( - \infty , 1)$ for $\displaystyle  i = 1, 2, ..., n$ and
 $$A(Y_t^x) =  p_0 {\bf 1}_{\{ Y_t^x \le a_1 \} } + \sum_{i=1}^{n-1} p_i {\bf 1}_{\{ a_i < Y_t^x \le a_{i+1} \} }  + p_N {\bf 1}_{\{ a_n < Y_t^x  \} } .$$
See \cite{BC} for the existence and uniqueness of strong solution of this SDE.
\end{remark}

\bigskip

{\bf Zhen-Qing Chen}

Department of Mathematics, University of Washington, Seattle,
WA 98195, USA

E-mail: \texttt{zqchen@uw.edu}

\bigskip

{\bf Mounir Zili}

Department of Mathematics, Faculty of sciences of Monastir, Tunisia

E-mail: \texttt{Mounir.Zili@fsm.rnu.tn}

\end{document}